\documentclass[a4paper, 12pt]{article}

\usepackage{amsfonts, amsmath, amsthm, amssymb}
\usepackage{latexsym} 

\theoremstyle{plain}
\newtheorem{thm}{Theorem}[section]
\newtheorem{prp}{Proposition}[section]
\newtheorem{lem}{Lemma}[section]
\newtheorem{cor}{Corollary}[section]

\theoremstyle{definition}
\newtheorem{dfn}{Definition}[section]

\theoremstyle{remark}
\newtheorem{rmk}{Remark}[section]

\makeatletter
 
 \@addtoreset{equation}{section}
\makeatother

\newcommand{\Z}{\mathbb{Z}}

\newcommand{\R}{\mathbb{R}}

\newcommand{\vomega}{\mbox{\boldmath$\omega$}}
\newcommand{\OOmega}{\mbox{\boldmath$\Omega$}}
\newcommand{\pa}{\partial}
\newcommand{\eps}{\varepsilon}

\newcommand{\nequiv}[1]{\overset{#1}{\simeq}}

\newcommand{\diag}{\rm diag}

\newcommand{\image}{{\rm Image}\,}

\newcommand{\fquad}{F^{\mbox{\tiny{\rm quad}}}}

\begin{document}
\title{
A note on the null condition for quadratic nonlinear 
  Klein-Gordon systems in two space dimensions
}  

\author{
         Soichiro Katayama\thanks{
             Department of Mathematics, Wakayama University. 
             930 Sakaedani, Wakayama 640-8510, Japan. 
             (E-mail: {\tt katayama@center.wakayama-u.ac.jp})} 
   \and  
         Tohru Ozawa\thanks{
             Department of Applied Physics, Waseda University. 
             3-4-1 Okubo, Shinjuku, Tokyo 169-8555, Japan. 
              (E-mail: {\tt txozawa@waseda.jp})}
   \and  
         Hideaki Sunagawa\thanks{
             Department of Mathematics, Osaka University. 
             1-1 Machikaneyama-cho, Toyonaka, Osaka 560-0043, Japan.
             (E-mail: {\tt sunagawa@math.sci.osaka-u.ac.jp})}
 } 
 
\date{ \today}   
\maketitle

\vspace{-4mm}
\begin{center}
Dedicated to the memory of Professor Yujiro Ohya\\
\end{center}
\vspace{2mm}

\begin{abstract}
 We consider the Cauchy problem for quadratic nonlinear 
 Klein-Gordon systems in two space dimensions with  masses 
 satisfying the resonance relation. 
 Under the null condition in the sense of 
 J.-M.~Delort, D.~Fang, R.~Xue (J.~Funct.~Anal.~{\bf 211} (2004), 288--323), 
 we show the global existence of asymptotically free solutions 
 if the initial data are sufficiently small in some weighted Sobolev space. 
 Our proof is based on an algebraic characterization of nonlinearities 
 satisfying the null condition. \\
 
\noindent{\bf Key Words:}\ 
Nonlinear Klein-Gordon systems; 
Null condition; Mass resonance. \\
\noindent{\bf Mathematics Subject Classification:}\ 
35L70, 35B40, 35L15\\
\noindent{\bf Running Title:}\ 
Null condition for  NLKG systems\\
\end{abstract}

\section{Introduction} 
In the present paper we consider large time behavior of solutions to 
the Cauchy problem for nonlinear systems of Klein-Gordon equations in two 
space dimensions:
\begin{align}
   (\Box + m_j^2) u_j = F_j(u,\pa u, \pa^2 u), 
 \qquad (t,x) \in [0,\infty)\times\R^{2}, \ j=1,2,
 \label{eq}
\end{align}
\begin{align}
  u_j(0,x)=f_j(x),\ \ \pa_t u_j(0,x)= g_j(x), 
  \qquad  x \in \R^2,\ j=1,2,
\label{data}
\end{align}
where $\Box=\pa_t^2-\pa_1^2-\pa_2^2$, $\pa=(\pa_0,\pa_1,\pa_2)$ with 
$\pa_0=\pa_t=\pa/\pa t$, $\pa_{j}=\pa/ \pa x_j$ for $j=1,2$, 
and $u=(u_k)_{k=1,2}$ is an $\R^2$-valued unknown function, 
while  
$\pa u=(\pa_a u_k)_{{\stackrel{\scriptstyle k=1,2}{\scriptstyle a=0,1,2}}}$ 
and 
$\pa^2 u=
(\pa_a \pa_b u_k)_{{\stackrel{\scriptstyle k=1,2}{\scriptstyle a,b=0,1,2}}}$ 
are its first and second order derivatives, respectively. 
The masses $m_1$, $m_2$ are positive constants. Without loss of generality 
we may always assume $m_1\le m_2$ throughout this paper. 
The nonlinear term $F_j =F_j(\xi, \eta, \zeta)$ 
is a $C^{\infty}$ function of 
$(\xi,\eta,\zeta) \in \R^2\times \R^{2\times 3}\times \R^{2\times 9}$ which 
vanishes of quadratic order at the origin, that is, 
$$
 F_j(\xi, \eta, \zeta)=O((|\xi|+|\eta|+|\zeta|)^2) 
 \quad \mbox{as} \quad 
 (\xi,\eta,\zeta) \to (0,0,0).
$$
We always suppose that the system is quasi-linear. In other words we assume 
that
\begin{equation}
\label{QL}
F_j(u,\pa u, \pa^2 u)
=
\sum_{k=1}^2\sum_{a,b=0}^3 \gamma^{jk}_{ab}(u,\pa u)\pa_a\pa_b u_k
+
\widetilde{F}_j(u,\pa u)
\end{equation}
with some functions $\gamma_{ab}^{jk}(\xi,\eta)$ vanishing of first order at 
the origin, and $\widetilde{F}_j(\xi,\eta)$ vanishing of quadratic order. 
To ensure the hyperbolicity, we assume that
\begin{equation}
\label{Symmetry}
\gamma_{ab}^{jk}(\xi, \eta)=\gamma_{ab}^{kj}(\xi,\eta),
\quad j,k=1,2,\ a,b=0,1,2,\ (\xi,\eta)\in \R^2\times \R^{2\times 3}.
\end{equation}
Without loss of generality, 
we may also assume that 
$\gamma_{a b}^{jk}(\xi,\eta)=\gamma_{ba}^{jk}(\xi,\eta)$ 
and  $\gamma_{00}^{jk}(\xi,\eta)\equiv 0$.

From the perturbative viewpoint, quadratic nonlinear interaction is of 
special interest for the Klein-Gordon equations in two space dimensions 
because large time  behavior of the solution is actually affected by 
the structure of the nonlinearities and by the ratio of the masses 
even if the initial data are sufficiently small, smooth and localized. 
In the case of $m_2 \ne 2m_1$ (which we call the non-resonant case), 
it is shown by Sunagawa~\cite{suna1} and Tsutsumi~\cite{tsutsumi} that 
the solution exists globally without any 
restrictions on $F_1$, $F_2$ (other than the hyperbolicity assumption) 
if $f_j$, $g_j$ are sufficiently small in 
a suitable weighted Sobolev space. Moreover, the solution is asymptotically 
free in the sense that we can find a solution $u^+(t)$ of the homogeneous 
linear Klein-Gordon equations satisfying
$$ 
 \lim_{t \to +\infty} \|u(t)-u^+(t)\|_{E}=0,
$$
where the energy norm $\|\cdot\|_E$ is defined by 
$$
 \|\phi(t)\|_{E}
 =
 \left(\sum_{j=1}^{2}\frac{1}{2} \int_{\R^2}
  (\pa_t \phi_j(t,x))^2 +|\nabla \phi_j(t,x)|^2+m_j^2\phi_j(t,x)^2\, dx
 \right)^{1/2}
$$
for $\phi=(\phi_j)_{j=1,2}$. 
On the other hand, in the resonant case (i.e., the case where $m_2=2m_1$) 
we must put some structural condition on the nonlinearities 
in order to obtain asymptotically free solution because there are 
examples of ($F_1$, $F_2$) such that the energy of the corresponding solution 
grows up as $t \to \infty$ (see \cite{suna2}). 
A sufficient condition on the nonlinearities is introduced by 
Delort-Fang-Xue \cite{dfx}, called the {\em null condition} (see  
Definition \ref{dfn_null} below), which allows us to show the global existence 
of small amplitude solutions for (\ref{eq})--(\ref{data}) in the resonant case 
if the data are sufficiently small, smooth and compactly supported
(see also Kawahara-Sunagawa~\cite{kawasuna}). 
A pointwise asymptotic profile of the global solution is also given in 
\cite{dfx}. However, their result does not imply the existence of a free 
profile in the  sense of the energy norm. 

The aim of this paper is to show the existence of a free profile 
in the sense of the energy norm under the null condition. 
Our approach is based on an algebraic characterization of the nonlinearities 
satisfying the null condition (see Proposition \ref{prp_characterize} below), 
which will allow us to reduce the problem essentially to 
the case of cubic nonlinearity
through a kind of normal form argument. Note also that, differently from 
\cite{dfx} and \cite{kawasuna}, our proof does not 
require compactness of the support of the initial data because we do not use 
the hyperbolic coordinates at all. We have only to assume that the initial 
data belong to some weighted Sobolev space and are sufficiently small in its 
norm, as in the non-resonant case \cite{suna1}. This is another advantage 
of our approach.

There is a large literature on the nonlinear Klein-Gordon equation. We refer 
the readers to \cite{delort}--\cite{geo2}, \cite{hn}--\cite{tsutsumi} and 
references therein.

\section{Main result} 

Let us first recall the definition of the null condition for the resonant 
quadratic nonlinear Klein-Gordon systems. We will follow the reformulation by 
Kawahara-Sunagawa (see the condition (a) in \cite{kawasuna}) instead of the 
original definition given by \cite{dfx}. 
For $j=1,2$, we denote by $\fquad_j$ the quadratic homogeneous part 
of $F_j$, 
that is, 
$$
 \fquad_j(\xi,\eta,\zeta)
 =\lim_{\lambda \downarrow 0} 
  \lambda^{-2} F_j(\lambda \xi, \lambda \eta, \lambda \zeta) 
$$
for $(\xi, \eta, \zeta) \in \R^2\times \R^{2\times 3}\times \R^{2\times 9}$. 
Also we set the unit hyperboloid
\begin{align*}
 \mathbb{H} =\{\vomega=(\omega_0,\omega_1,\omega_2) \in \R^3\, :\, 
  \omega_0^2-\omega_1^2-\omega_2^2=1  \}
\end{align*}
and 
\begin{equation}
 \Phi_j(\vomega)
 =\int_{0}^{1} 
  \fquad_j \bigl( U(\theta), V(\vomega,\theta), W(\vomega,\theta) \bigr)\, 
 e^{-2 \pi i j\theta}\, d\theta
\label{NullKG}
\end{equation}
for $\vomega \in \mathbb{H}$, where 
\begin{align}
 &U(\theta) = \bigl(\cos 2 \pi k\theta \bigr)_{k=1,2}, 
 \label{U}\\ 
 &V(\vomega,\theta) = 
  \bigl( -\omega_a m_k \sin 2\pi k\theta 
  \bigr)_{{\stackrel{\scriptstyle k=1,2}{\scriptstyle a=0,1,2}}},
  \label{V}\\ 
 &W(\vomega, \theta)
  =\bigl(
   -\omega_a\omega_b m_k^2 \cos 2\pi k\theta
 \bigr)_{{\stackrel{\scriptstyle k=1,2}{\scriptstyle a,b=0,1,2}}}
 \label{W}
\end{align} 
and $i=\sqrt{-1}$. 

\begin{dfn}\label{dfn_null}
We say that the nonlinear term $(F_1,F_2)$ satisfies the null condition if 
$\Phi_1(\vomega)=\Phi_2(\vomega)=0$ for all  $\vomega \in \mathbb{H}$. 
\end{dfn}
Examples of quasi-linear term $(F_1,F_2)$ which satisfies the null condition 
as well as the hyperbolicity assumption \eqref{Symmetry} will be given in 
Remark \ref{rmk_example} below.

In order to state the main result precisely, let us also introduce 
the weighted Sobolev space as follows: 
$$
 H^{s,k}(\R^2)
 =
 \{ 
  \phi \in L^2(\R^2)\, :\, 
  (1+|\cdot|^2)^{k/2} (1-\Delta)^{s/2} \phi \in L^2(\R^2) 
 \}
$$
equipped with the norm 
$$
 \|\phi\|_{H^{s,k}(\R^2)}
 = \|(1+|\cdot|^2)^{k/2} (1-\Delta)^{s/2} \phi\|_{L^2(\R^2)}.
$$
For simplicity, 
we write $H^s=H^{s,0}$ and $\|\phi\|_{H^s}=\|\phi\|_{H^{s,0}}$. 

Now we are in a position to state the main result.
\begin{thm}\label{thm_main}
Let $m_2=2m_1$. Assume that $(F_1,F_2)$ satisfies the null condition in the 
sense of Definition $\ref{dfn_null}$. Assume also that 
\eqref{QL}--\eqref{Symmetry} is satisfied. 
Let $f_j\in H^{s+1,s}(\R^2)$, $g_j \in H^{s,s}(\R^2)$ for $j=1,2$ with 
$s \ge 29$. There exists a positive constant $\eps$ such that if 
\begin{equation}
\label{DataSize}
 \sum_{j=1}^{2}
 \Bigl(\|f_j\|_{H^{s+1,s}(\R^2)}+ \|g_j\|_{H^{s,s}(\R^2)}\Bigr) \le \eps,
\end{equation}
the problem \eqref{eq}--\eqref{data} admits a unique global solution 
$u=(u_1,u_2)$ satisfying 
$$
 u\in \bigcap_{k=0}^{s}C^{k}([0,\infty);H^{s+1-k}(\R^2)).
$$
Furthermore, $u$ has a free profile, i.e., there exists 
$(f_j^{+},g_j^{+}) \in H^{s-3}(\R^2)\times H^{s-4}(\R^2)$ 
such that 
$$
\lim_{t\to + \infty}
\Bigl(\|(u_j-u_j^{+})(t,\cdot)\|_{H^{s-3}(\R^2)} 
+ \|\pa_t(u_j-u_j^{+})(t,\cdot)\|_{H^{s-4}(\R^2)} \Bigr)
=0,
$$
where $u_j^{+}$ solves 
$(\Box+m_j^2)u_j^{+}=0$ 
with $(u_j^{+},\pa_{t}u_j^{+})|_{t=0}=(f_j^{+},g_j^{+})$ 
for $j=1,2$.\\
\end{thm}

\begin{rmk}
In the paper by Kawahara-Sunagawa \cite{kawasuna}, another sufficient 
condition for global existence with $C_0^{\infty}$ small data 
is also introduced (see the condition (b) in \cite{kawasuna}). 
Our proof in the present paper does not work for that case, 
and as pointed out in \cite{kawasuna}, 
some long-range effect should be taken into account. 
It is still an open problem to find out precise asymptotic profile of 
the global solution under the condition (b) even in the simplest case 
$$
 \left\{
 \begin{array}{l}
  F_1=u_1u_2,\\ F_2 =u_1^2, 
 \end{array}
 \right.
$$
that is a typical example which satisfies the 
condition (b) but violates the null condition. 
For closely related works on nonlinear Schr\"odinger systems, 
see the recent papers by Hayashi-Li-Naumkin \cite{hln}, \cite{hln2} 
and by Hayashi-Li-Ozawa \cite{hlo}.  \\
\end{rmk}

The rest of this paper is organized as follows: 
In the next section, we give some preliminaries mainly on the commuting vector 
fields and the null forms. In Section 4, we recall an algebraic normal form 
transformation developed in the previous papers. A characterization of the 
nonlinearities satisfying the null condition will be given in Section 5. 
After that, we will prove the main theorem in Section 6. 
Throughout this paper, we will frequently use the following 
conventions on implicit constants: 
\begin{itemize}
\item
$A\lesssim B$ (resp. $A \gtrsim B$) stands for  $A \le CB$ (resp. $A \ge CB$) 
with a positive constant $C$. 
\item 
The expression $f=\sum_{a \in A}' g_a$ means that there exists a family 
$\{C_a\}_{a \in A}$ of real constants such that $f=\sum_{a \in A} C_a g_a$. 
\end{itemize}
Also, the notation $\langle y \rangle=(1+|y|^2)^{1/2}$ will be used for 
$y \in \R^N$ with a positive integer $N$. 

\section{Notations and preliminaries}

We put
$x_0=-t$, $x=(x_1,x_2)$,  
$\Omega_{ab} = x_{a}\pa_{b}-x_{b}\pa_{a}$, $0\le a, b \le 2$, and 
$$
 Z=(Z_1,\dots, Z_{6})
 =\bigl( \pa_0,\pa_1,\pa_2,\Omega_{01}, \Omega_{02}, \Omega_{12} \bigr).
$$
Note that the following commutation relations hold: 
\begin{align}
&[\Box +m^2, Z_j] = 0, 
\label{commute}\\
&[\Omega_{a b}, \pa_{c}]  = \eta_{b c} \pa_{a} - \eta_{c a} \pa_{b}, 
\nonumber\\
&[\Omega_{a b}, \Omega_{c d}] 
  =  \eta_{a d} \Omega_{b c}+ \eta_{b c} \Omega_{a d}
   - \eta_{a c} \Omega_{b d} - \eta_{b d} \Omega_{a c}
\nonumber
\end{align}
for $m \in \R$, $1\leq j\leq 6$, $0\leq a, b \leq 2$. 
Here 
$[\cdot ,\cdot ]$ denotes the commutator of linear operators,
and 
$(\eta_{a b})_{0\leq a,b \leq 2}= \diag (-1,1,1)$. 
For a smooth function $\phi$ of $(t,x) \in \R^{1+2}$ and for a non-negative 
integer $s$, we define 
$$
|\phi(t,x)|_{s} := 
   \sum_{|\alpha|\leq s}|Z^{\alpha}\phi(t,x)|
$$
and 
$$
\|\phi(t)\|_{s} := 
   \sum_{|\alpha|\leq s}\|Z^{\alpha} \phi(t,\cdot)\|_{L^2(\R^2)}, 
$$
where $\alpha=(\alpha_1,\dots, \alpha_{6})$ is a multi-index, 
$Z^{\alpha} =Z_1^{\alpha_1}\cdots Z_{6}^{\alpha_{6}}$ and 
$|\alpha|= \alpha_{1}+\cdots +\alpha_{6}$. 
Next we introduce the null forms 
\begin{align}
  &Q_{a b}(\phi,\psi)
 =(\pa_{a} \phi) (\pa_{b} \psi)-(\pa_{b} \phi) (\pa_{a} \psi),
 \quad  0\le a,b\le 2,
  \label{strongnull}\\ 
 &Q_0(\phi,\psi)
 =(\pa_t \phi)(\pa_t \psi) -(\nabla_x \phi )\cdot (\nabla_x \psi)
 =-\sum_{a,b=0}^{2} \eta_{ab}(\pa_a \phi)(\pa_b \psi).
 \label{Yowanull}
\end{align}
As pointed out in \cite{geo} (see also \cite{katayama2}), 
$Q_{a b}$ has a certain compatibility with the Klein-Gordon operator 
$\Box +m^2$, while $Q_0$ is sometimes not so if $m\ne 0$. 
In what follows, we call $Q_{ab}$ {\em the strong null forms}.  
From the identities 
\begin{align*}
 {x}_{c}Q_{a b} (\phi,\psi) 
  &= (\Omega_{c a}\phi) (\pa_{b} \psi) + (\Omega_{b c}\phi) (\pa_{a} \psi)
     + (\pa_{c} \phi) (\Omega_{a b} \psi), \\
\pa_{c}Q_{a b} (\phi,\psi) 
  &= Q_{a b}(\pa_{c}\phi,\psi) + Q_{a b}(\phi,\pa_{c}\psi), \\
\Omega_{c d} Q_{a b} (\phi,\psi) 
  &= Q_{a b}(\Omega_{c d}\phi,\psi) 
     + Q_{a b}(\phi,\Omega_{c d}\psi)
     + \eta_{a c}Q_{b d} (\phi,\psi)\\
  &\quad  
     + \eta_{b d}Q_{a c}(\phi,\psi) 
     - \eta_{a d}Q_{b c}(\phi,\psi) 
     - \eta_{b c}Q_{a d}(\phi,\psi),
\end{align*}
we deduce the following properties on the strong null forms. 

\begin{lem} \label{lem_null}
Let $\phi$, $\psi$ be smooth functions on $(t,x) \in \R^{1+2}$. We have 
\begin{align*}
 |Q_{a b}(\phi,\psi)|  \lesssim \frac{1}{\langle t+|x| \rangle} 
   \bigl( |\phi|_1 |\pa \psi| + |\pa \phi| |\psi|_1 \bigr) 
\end{align*}
for $0\le a,b\le 2$, and
\begin{align*}
 Z^{\alpha} Q_{a b} (\phi,\psi) = 
  \sum_{c,d =0}^{2} \mathop{{\;\,\sum}'}_{|\beta|+|\gamma|\leq |\alpha|}
  Q_{c d} (Z^{\beta} \phi, Z^{\gamma} \psi) 
\end{align*}
for any multi-index $\alpha$.\\
\end{lem}

We close this section with the following decay estimate due to Georgiev
\cite{geo2}.

\begin{lem} \label{lem_decay}
Let $m$ be a positive constant and $w$ be a solution of 
the inhomogeneous linear Klein-Gordon equation 
$(\square +m^2)w=h$ for $t\geq 0$, $x \in \R^2$. 
Then we have
\begin{align*}
 \langle t+|x|\rangle |w(t,x)|
  \lesssim &  
     \sum_{j=0}^{\infty} \sum_{|\beta|\leq 4}
      \sup_{ \tau \in [0,t]} \varphi_j(\tau) 
     \bigl\| 
          \langle \tau +|y| \rangle Z^{\beta} h(\tau ,y ) 
     \bigr\| _{L^2(\R_y^2)}\\
  &+\sum_{j=0}^{\infty} \sum_{|\beta|\leq 5} 
      \left\| 
        \langle y \rangle \varphi_j(|y|) \bigl(Z^{\beta} w \bigr) (0,y )
      \right\|_{L^2(\R_y^2)},
\end{align*}
provided that the right-hand side is finite. 
Here $\{\varphi_j\}_{j=0}^{\infty}$ is the Littlewood-Paley partition of 
unity, i.e.,
\begin{align*}
  \sum_{j=0}^{\infty} \varphi_j(\tau) =1 \ \  (\tau \geq 0); \quad 
  \varphi_j \in C_0^{\infty}(\R),\ 
  \varphi_j \geq 0 \ \mbox{for}\ j\geq 0; \\
  {\rm supp}\ \varphi_j \subset [2^{j-1},2^{j+1}] \mbox{ for } j\geq 1, \ \ 
  {\rm supp}\ \varphi_0 \cap \R_{\ge 0} \subset [0,2].\\
\end{align*}
\end{lem}

\section{Algebraic normal form transformation}

In this section, we recall an algebraic normal form transformation developed 
by \cite{kosecki}, \cite{katayama}, \cite{tsutsumi}, \cite{suna1}, etc. 
Let $v_j$ and $\widetilde{v}_j$ be smooth functions on $(t,x) \in \R^{1+2}$. 
We set 
$h_j=(\Box+m_j^2)v_j$ and $\widetilde{h}_j=(\Box+m_j^2)\widetilde{v}_j$ 
for $j=1,2$. 
Throughout this section, we use the the following convention: 
We write 
$$\Phi\sim \Psi$$ 
if $\Phi-\Psi$ can be written as a linear combination of 
$Q_{ab}(\pa^\alpha v_k, \pa^\beta \widetilde{v}_l)$, 
$(\pa^\alpha v_k) (\pa^{\beta} \widetilde{h}_l)$, 
$(\pa^\alpha h_k)(\pa^\beta\widetilde{v}_l)$ or $h_k\widetilde{h}_l$ 
with $|\alpha|$, $|\beta| \le 1$, $0\le a,b \le 2$ and $1\le k, l \le 2$, 
where $Q_{ab}$ is given by \eqref{strongnull}.

\begin{prp} \label{prp_algebra}
Put $\mathbf{e}_{kl}=v_k \widetilde{v}_l$, 
$\mathbf{\tilde{e}}_{kl}=Q_0(v_k, \widetilde{v}_l)$ 
and $\mathcal{L}_j=\Box+m_j^2$, where $Q_0$ is given by \eqref{Yowanull}. 
We have 
$$
 (\mathcal{L}_j\mathbf{e}_{kl}\ \ \mathcal{L}_j\mathbf{\tilde{e}}_{kl}) 
 \sim 
 (\mathbf{e}_{kl}\ \ \mathbf{\tilde{e}}_{kl}) A_{jkl},
$$
where 
\begin{align*} 
A_{jkl}=\begin{pmatrix}
 m_j^2-m_k^2-m_l^2 & 2m_k^2m_l^2\\
 2 & m_j^2-m_k^2-m_l^2 
 \end{pmatrix}.
\end{align*}
\end{prp}

\noindent {\bf Proof.}\ \ 
This proposition is nothing but a paraphrase of Lemma 6.1 of \cite{suna1}. 
However, for the convenience of the readers, we give a proof here. 
It is sufficient to show that  
\begin{align} 
 \Box \mathbf{e}_{kl}
 \sim -(m_k^2+m_l^2) \mathbf{e}_{kl}+2\mathbf{\tilde{e}}_{kl}
\label{first}
\end{align}
and
\begin{align} 
 \Box \mathbf{\tilde{e}}_{kl}
 \sim 
 -(m_k^2+m_l^2)\mathbf{\tilde{e}}_{kl}+2m_k^2m_l^2 \mathbf{e}_{kl}.
\label{second}
\end{align}
Since $\Box v_k=-m^2_k v_k + h_k$ and 
$\Box \widetilde{v}_l=-m_l^2 \widetilde{v}_l+\widetilde{h}_l$, we have 
\begin{align*}
 \Box \mathbf{e}_{kl}
 &=
 (\Box v_k)\widetilde{v}_l + 
 v_k \Box \widetilde{v}_l+2 Q_0(v_k, \widetilde{v}_l)\\
 &\sim
 (-m_k^2 v_k)\widetilde{v}_l+
 v_k(-m_l^2 \widetilde{v}_l) +2 \mathbf{\tilde{e}}_{kl}\\
 &=
 -(m_k^2+m_l^2) \mathbf{e}_{kl}+ 2 \mathbf{\tilde{e}}_{kl},
\end{align*}
which yields the first relation (\ref{first}). 
As for the second relation (\ref{second}), 
we observe the relation
$$
 (\pa_a\pa_c v_k)(\pa_b \pa_d \widetilde{v}_l)
 =
 (\pa_a\pa_b v_k)(\pa_c \pa_d \widetilde{v}_l) + Q_{cb}(\pa_a v_k, \pa_d 
\widetilde{v}_l)
 \sim 
 (\pa_a\pa_b v_k)(\pa_c \pa_d \widetilde{v}_l)
$$
to obtain 
\begin{align*}
 \Box \mathbf{\tilde{e}}_{kl}
 &=
 Q_0(\Box v_k, \widetilde{v}_l) + Q_0(v_k, \Box \widetilde{v}_l) 
 +
2\sum_{a,b} \sum_{c,d} 
  \eta_{ab} \eta_{cd} (\pa_a\pa_{c} v_k)(\pa_b \pa_d \widetilde{v}_l)\\
 &\sim
 Q_0(-m_k^2 v_k, \widetilde{v}_l) +Q_0(v_k,-m_l^2 \widetilde{v}_l) 
 +
2\sum_{a,b} \sum_{c,d} \eta_{ab} 
  \eta_{cd} (\pa_a\pa_{b}v_k)(\pa_c \pa_d \widetilde{v}_l)\\
 &=
 -(m_k^2+m_l^2) Q_0(v_k, \widetilde{v}_l)+ 2 (\Box v_k)(\Box \widetilde{v}_l)\\
 &\sim
 -(m_k^2+m_l^2) \mathbf{\tilde{e}}_{kl} 
 + 2 (-m_k^2 v_k)(-m_l^2 \widetilde{v}_l)\\
 &=
 -(m_k^2+m_l^2) \mathbf{\tilde{e}}_{kl} + 2m_k^2 m_l^2 \mathbf{e}_{kl}.
\end{align*}
This completes the proof. \qed\\
\medskip

Now we focus our attention on the structure of the matrix $A_{jkl}$ for 
$1\le j\le 2$ and $1\le k \le l \le 2$ under the resonance relation 
$m_2=2m_1>0$. Since 
$$
 \det A_{jkl}
 =
 \prod_{\sigma_1,\sigma_2 \in \{\pm 1\}} (m_j+\sigma_1 m_k+\sigma_2 m_l),
$$
we can see that $A_{jkl}$ is invertible if and only if 
$(j,k,l) = (1,1,1)$, $(1,2,2)$, $(2,1,2)$ or $(2,2,2)$. In this case, 
we have 
\begin{align}
 v_k \widetilde{v}_l
 &=
 (\mathbf{e}_{kl}\ \ \mathbf{\tilde{e}}_{kl})
 \begin{pmatrix}1 \\ 0 \end{pmatrix} 
 \nonumber\\
 &=
 (\mathbf{e}_{kl}\ \ \mathbf{\tilde{e}}_{kl})
 A_{jkl}\begin{pmatrix}p_{jkl} \\ \tilde{p}_{jkl} \end{pmatrix}
 \nonumber\\
 &\sim
  (\mathcal{L}_j\mathbf{e}_{kl}\ \ \mathcal{L}_j\mathbf{\tilde{e}}_{kl})
 \begin{pmatrix}p_{jkl} \\ \tilde{p}_{jkl} \end{pmatrix}
 \nonumber\\
  &=
  (\Box+m_j^2)\Bigl(p_{jkl} v_k\widetilde{v}_l 
  + \tilde{p}_{jkl} Q_0(v_k, \widetilde{v}_l) \Bigr)
\label{NormalF01}
\end{align}
with 
$$
 \begin{pmatrix} p_{jkl}\\ \tilde{p}_{jkl} \end{pmatrix}
 =A_{jkl}^{-1} \begin{pmatrix} 1\\ 0 \end{pmatrix}.
$$
On the other hand, $A_{jkl}$ is degenerate when $(j,k,l) = (1,1,2)$ or 
$(2,1,1)$. Indeed, 
$$
 \image (A_{112})=\left\{ 
 \kappa \begin{pmatrix} -2m_1^2\\ 1
 \end{pmatrix}
 \ :\ 
 \kappa \in \R
 \right\},
 \quad 
 \image (A_{211})=\left\{ 
 \kappa \begin{pmatrix} m_1^2\\ 1
 \end{pmatrix}
 \ :\ 
 \kappa \in \R
 \right\}.
$$
Remark that $\image (A_{112})$ and $\image (A_{211})$ correspond directly to 
the quadratic terms 
\begin{align}
 G_1(v_1, \widetilde{v}_2)
 :=
 Q_0(v_1, \widetilde{v}_2)-2m_1^2 v_1 \widetilde{v}_2
 \sim
 (\Box +m_1^2)\left(\frac{v_1 \widetilde{v}_2}{2}\right)
 \label{g1}
\end{align}
and
\begin{align}
 G_2(v_1, \widetilde{v}_1)
 :=
 Q_0(v_1, \widetilde{v}_1) + m_1^2 v_1\widetilde{v}_1 
 \sim
  (\Box +m_2^2)\left(\frac{v_1\widetilde{v}_1}{2}\right),
  \label{g2}
\end{align}
respectively.
The following quadratic terms should be also associated with 
$\image (A_{112})$: 
\begin{align}
 H_{1,a}(v_1, \widetilde{v}_2):=v_1 \pa_a \widetilde{v}_2 + 
2\widetilde{v}_2\pa_a v_1, \quad a=0,1,2.
  \label{h1a}
\end{align}
In fact, it follows that 
\begin{align*}
2m_1^2H_{1,a}(v_1,\widetilde{v}_2)+G_1(v_1, \pa_a \widetilde{v}_2)
=& Q_0(v_1, \pa_a\widetilde{v}_2)+4m_1^2(\pa_a v_1)\widetilde{v}_2\\
=& - \sum_{b,c}\eta_{bc}(\pa_b v_1)(\pa_c\pa_a \widetilde{v}_2)
   + m_2^2(\pa_a v_1)\widetilde{v}_2\\ 
\sim & 
  -\sum_{b,c} \eta_{bc}(\pa_a v_1)(\pa_b\pa_c\widetilde{v}_2)
  + m_2^2(\pa_a v_1)\widetilde{v}_2\\
= & (\pa_a v_1)\widetilde{h}_2\sim 0,
\end{align*}
whence 
\begin{equation}
 H_{1, a}(v_1, \widetilde{v}_2) 
 \sim 
 -\frac{1}{2m_1^2} G_1(v_1, \pa_a \widetilde{v}_2)
 \sim  
 (\Box +m_1^2)\left(\frac{v_1\pa_a \widetilde{v}_2}{-4m_1^2}\right).
\label{NormalF02}
\end{equation}
Similarly
\begin{equation}
\label{h2a}
 H_{2,a}(v_1,\widetilde{v}_1)
 :=v_1 \pa_a \widetilde{v}_1 - \widetilde{v}_1 \pa_a v_1, \quad a=0,1,2, 
\end{equation}
are associated with $\image (A_{211})$ since we can see that 
\begin{equation}
H_{2,a}(v_1,\widetilde{v}_1)
\sim 
\frac{1}{m_1^2} G_2(v_1, \pa_a \widetilde{v}_1)
 \sim  
 (\Box +m_2^2)\left(\frac{v_1\pa_a \widetilde{v}_1}{2m_1^2}\right).
\label{NormalF03}
\end{equation}
The above observation will play the key role in the proof of 
Theorem \ref{thm_main}.

\section{Characterization of the null condition}

The aim of this section is to give a characterization of the null condition 
in terms of $G_1$, $G_2$, $H_{1, a}$, $H_{2,a}$, and $Q_{ab}$ defined 
in the previous sections. What we are going to prove is the following. 

\begin{prp} \label{prp_characterize}
$(F_1, F_2)$ satisfies the null condition if and only if its quadratic 
homogeneous part can be  written in the following form: 
\begin{align}
 \fquad_1(u,\pa u, \pa^2 u)
 =&
 \mathop{{\;\,\sum}'}_{|\alpha|+|\beta|\le 1}
  G_1(\pa^\alpha u_1, \pa^\beta u_2)
+ \sum_{a=0}^2 \mathop{{\;\,\sum}'}_{|\alpha|, |\beta|\le 1} 
  H_{1,a}(\pa^\alpha u_1, \pa^\beta u_2)
\nonumber\\
&+\sum_{a,b=0}^2\mathop{{\;\,\sum}'}_{|\alpha|+|\beta|\le 1} 
 Q_{ab}(\pa^\alpha u_1, \pa^\beta u_2) \nonumber\\
&+ 
 \mathop{{\;\,\sum}'}_{|\alpha|, |\beta|\le 2;\, |\alpha|+|\beta|\le 3} 
 (\pa^{\alpha} u_1)(\pa^{\beta} u_1) 
 + 
 \mathop{{\;\,\sum}'}_{|\alpha|, |\beta|\le 2;\, |\alpha|+|\beta|\le 3} 
 (\pa^{\alpha} u_2)(\pa^{\beta} u_2), 
 \label{f1null}
\end{align}
\begin{align}
 \fquad_2(u,\pa u, \pa^2 u)
 = &
 \mathop{{\;\,\sum}'}_{|\alpha|\le 1} G_2(u_1,\pa^\alpha u_1) 
 +
 \mathop{{\;\,\sum}'}_{a,b=0}^2 H_{2,a}(u_1, \pa_b u_1)
 +\mathop{{\;\,\sum}'}_{a,b,c=0}^2 Q_{ab}(u_1, \pa_c u_1)\nonumber\\
 &+\mathop{{\;\,\sum}'}_{|\alpha|, |\beta|\le 2;\, |\alpha|+|\beta|\le 3} 
 (\pa^{\alpha} u_1)(\pa^{\beta} u_2)
 +\mathop{{\;\,\sum}'}_{|\alpha|, |\beta|\le 2;\, |\alpha|+|\beta|\le 3} 
 (\pa^{\alpha} u_2)(\pa^{\beta} u_2),
 \label{f2null}
\end{align}
where $G_1$, $G_2$, $H_{1,a}$, $H_{2,a}$, and $Q_{ab}$ 
are given by \eqref{g1}, \eqref{g2}, \eqref{h1a}, \eqref{h2a}, and 
\eqref{strongnull}, respectively.
\end{prp}

\begin{rmk} \label{rmk_example}
If we further assume \eqref{QL}--\eqref{Symmetry} in addition to the null 
condition, some restriction on the coefficients for the terms including 
$\pa^2 u$ is needed. 
For example, the following $(F_1, F_2)$ satisfies the null condition 
as well as \eqref{QL}--\eqref{Symmetry}:
\begin{align*}
F_1
=
&\sum_{a=1}^2 p_{a} G_1(u_1, \pa_a u_2)
 + 
 \sum_{0\le a,b \le 2; a+b \ne 0} q_{ab} H_{1,a}(u_1, \pa_b u_2)\\
&{}
  + \mathop{{\;\,\sum}'}_{a=1}^2  G_1(\pa_a u_1, u_2) 
  + \mathop{{\;\,\sum}'}_{0\le a,b \le 2; a+b \ne 0} H_{1,a}(\pa_b u_1,u_2),\\
F_2
=
&\sum_{a=1}^2 p_{a}G_2(u_1, \pa_a u_1)
 + \sum_{0\le a,b \le 2; a+b \ne 0} q_{ab} H_{2,a}(u_1, \pa_b u_1)
\end{align*}
with real constants $p_{a}$ and $q_{ab}$.
\end{rmk}

\noindent{\bf Proof of Proposition~\ref{prp_characterize}.}\ \ 
For $j=1$, $2$, we will write 
$$
 \Psi(u,\pa u, \pa^2 u) \nequiv{j}\widetilde{\Psi}(u,\pa u, \pa^2 u)
$$ 
 if we have
$$
\int_0^1 
\bigl(
 \Psi-\widetilde\Psi\bigr)
 \bigl(U(\theta), V(\vomega,\theta), W(\vomega,\theta)
\bigr) 
 e^{-2\pi i j\theta} d\theta=0
$$
for $\vomega\in \mathbb{H}$, where $U$, $V$ and $W$ are given by \eqref{U}, 
\eqref{V} and \eqref{W}, respectively. 
We split $\fquad_j(u,\pa u, \pa^2 u)$ into the three parts: 
$$
\fquad_j(u,\pa u, \pa^2 u)
=
F_j^{11}(u,\pa u,\pa^2 u)+F_j^{12}(u,\pa u, \pa^2 u)+F_j^{22}(u,\pa u, \pa^2 u)
$$
with
$$
F_j^{kl}(u,\pa u, \pa^2 u)
=
\mathop{{\;\,\sum}'}_{|\alpha|, |\beta|\le 2;\, |\alpha|+|\beta|\le 3} 
(\pa^\alpha u_k)(\pa^\beta u_l).
$$
Then we can check that
$$
\fquad_1(u,\pa u, \pa^2 u) \nequiv{1} F_1^{12}(u, \pa u, \pa^2 u),
\quad \fquad_2(u,\pa u, \pa^2 u) \nequiv{2} F_2^{11}(u, \pa u, \pa^2 u)
$$
by using the relation 
$$
\int_0^1 e^{2\pi i (k_1+k_2-j)\theta} d\theta
=\left\{
 \begin{array}{ll}
 1 & \mbox{if \  $(j,k_1,k_2)=(1,-1,2)$\ or\ $(2,1,1)$},\\
 0 & \mbox{otherwise},
 \end{array}
\right.
$$
for $j=1,2$ and $k_1,k_2 \in \Z$ with $1\le |k_1|\le |k_2| \le 2$. 
Hence we only have to investigate $F_1^{12}$ and $F_2^{11}$ in order to check 
the null condition.

First we consider $F_1^{12}$. We rewrite $(\pa^\alpha u_1)(\pa^\beta u_2)$ 
with $|\alpha|=0$ as follows:
\begin{align*}
 u_1(\pa^\beta u_2)
 =& 
 \frac{1}{2m_1^2} \left(Q_0(u_1, \pa^\beta u_2)-G_1(u_1, \pa^\beta u_2)\right),
 \quad |\beta|\le 1, \\ 
 u_1(\pa_a\pa_b u_2)
 = & 
 -2(\pa_a u_1)(\pa_b u_2)+H_{1,a}(u_1, \pa_b u_2).
\end{align*}
We can rewrite $(\pa^\alpha u_1)(\pa^\beta u_2)$ with $|\beta|=0$ 
in a similar fashion. Moreover we have 
$$
 (\pa_c u_1)(\pa_a\pa_b u_2)
 =
 -2(\pa_b u_2)(\pa_a\pa_c u_1) + H_{1,a}(\pa_cu_1, \pa_b u_2).
$$
Hence we find
\begin{align*}
F_1^{12}(u,\pa u, \pa^2 u)
=& 
 \sum_{a,b=0}^2 \lambda_{ab} (\pa_a u_1)(\pa_b u_2)
 +\sum_{a,b,c=0}^2 \mu_{abc}(\pa_c u_2)(\pa_a\pa_b u_1)\\
& + \mathop{{\;\,\sum}'}_{|\alpha|+|\beta|\le 1} 
     G_1(\pa^\alpha u_1, \pa^\beta u_2)
  + \sum_{a=0}^2\mathop{{\;\,\sum}'}_{|\alpha|,|\beta|\le 1} 
     H_{1,a}(\pa^\alpha u_1, \pa^\beta u_2)
\end{align*}
with some real constants $\lambda_{ab}$ and $\mu_{abc}$. 
Since $G_1(\pa^\alpha u_1, \pa^\beta u_2)\nequiv{1} 0$ 
and $H_{1,a}(\pa^\alpha u_1, \pa^\beta u_2) \nequiv{1} 0$, 
it follows from the definition of $\Phi_1$ that
\begin{align*}
 \Phi_1(\vomega)
 =\frac{m_1m_2}{4} \sum_{a,b=0}^2 \lambda_{ab}\omega_a\omega_b
 -i\frac{m_1^2m_2}{4} \sum_{a,b,c=0}^3 \mu_{abc} \omega_a\omega_b\omega_c.
\end{align*}
In order that this quantity vanishes identically on $\mathbb{H}$, 
we must have $\lambda_{aa}=0$ for $0\le a\le 2$, 
$\lambda_{ab}=-\lambda_{ba}$ for $0\le a<b\le 2$, and so on.
Now it is not difficult to see that
$\sum_{a,b} \lambda_{ab} (\pa_a u_1)(\pa_b u_2)$
and $\sum_{a,b,c}\mu_{abc} (\pa_c u_2)(\pa_a\pa_b u_1)$
can be written in terms of the strong null forms.
Hence we have (\ref{f1null}). The converse is also true. 
Similarly, by writing 
\begin{align*}
u_1(\pa^\alpha u_1)
=& 
-\frac{1}{m_1^2}\left(Q_0(u_1,\pa^\alpha u_1)-G_2(u_1, \pa^\alpha u_1)\right), 
\quad |\alpha|\le 1,\\
u_1(\pa_a\pa_b u_1)
=& 
(\pa_a u_1)(\pa_b u_1) + H_{2,a}(u_1, \pa_b u_1),
\end{align*}
we have
\begin{align*}
F_2^{11}(u,\pa u, \pa^2 u)
=& 
 \sum_{a,b=0}^2 \tilde\lambda_{ab} (\pa_a u_1)(\pa_b u_1)
 +\sum_{a,b,c=0}^2 \tilde\mu_{abc}(\pa_c u_1)(\pa_a\pa_b u_1)\\
& + \mathop{{\;\,\sum}'}_{|\alpha|\le 1} G_2(u_1, \pa^\alpha u_1)
  +\mathop{{\;\,\sum}'}_{a, b=0}^2 H_{2,a}(u_1, \pa_b u_1)
\end{align*}
with appropriate real constants $\tilde{\lambda}_{ab}$ and $\tilde\mu_{abc}$, 
which leads to
\begin{align*}
 \Phi_2(\vomega)
 =&
 -\frac{m_1^2}{4}\sum_{a,b=0}^{2}\tilde{\lambda}_{ab} \omega_a \omega_b
 -i\frac{m_1^3}{4}\sum_{a,b,c=0}^2 \tilde{\mu}_{abc} \omega_a \omega_b\omega_c.
\end{align*}
As before, $\Phi_2(\vomega)\equiv 0$ on $\mathbb{H}$
implies that $\sum_{a,b}\tilde{\lambda}_{ab}(\pa_a u_1)(\pa_b u_1)$ and
$\sum_{a,b,c}\tilde\mu_{abc}(\pa_cu_1)(\pa_a\pa_b u_1)$ can be written 
in terms of the strong null forms. This leads to (\ref{f2null}). 
The converse is also true. \qed\\

By Proposition~\ref{prp_characterize}, \eqref{NormalF01}, \eqref{g1},
\eqref{g2}, \eqref{NormalF02}, and \eqref {NormalF03} we obtain the following.
\begin{cor} \label{cor_normalform}
Let $(u_1,u_2)$ be a smooth solution for \eqref{eq}. 
If $(F_1, F_2)$ satisfies the null condition, we have  
$$
 F_j(u,\pa u, \pa^2 u)  =  (\Box+m_j^2) \Lambda_j + N_j + R_j
$$
for $j=1,2$, where 
\begin{align*} 
 &\Lambda_j=
  \sum_{k,l=1}^{2}
  \mathop{{\;\,\sum}'}_{|\alpha|,|\beta|\le 3}
  (\pa^{\alpha}u_k)(\pa^{\beta} u_l),\\
 &N_j=
  \sum_{k,l=1}^{2}\sum_{a,b=0}^{2}
  \mathop{{\;\,\sum}'}_{|\alpha|,|\beta|\le 3}
  Q_{ab}(\pa^{\alpha}u_k,\pa^{\beta} u_l),
\end{align*}
and $R_j$ is a smooth function of $(\pa^\alpha u)_{|\alpha|\le 5}$ with
$$
R_j=
  O\Bigl( \bigl| (\pa^{\alpha} u)_{|\alpha|\le 5} \bigr|^3 \Bigr)
\quad \text{near $(\pa^\alpha u)_{|\alpha|\le 5}=0$}.
$$
\end{cor}

\section{Proof of the main theorem}

Now we are ready to prove Theorem~\ref{thm_main}. 
The main step of the proof is to get some {\it a priori} estimate. 
From now on, we suppose that the null condition, 
as well as \eqref{QL}--\eqref{Symmetry}, is satisfied and 
let $u=(u_j)_{j=1,2}$ be a solution to (\ref{eq})--(\ref{data}) 
for $t \in [0,T)$. We define 
\begin{align*}
 E(T)=
 \sup_{0\leq t < T} & 
  \Bigl[\ 
    \langle t\rangle^{-\delta}
      \bigl( \| u(t) \|_{s}+ \|\pa u(t)\|_{s}\bigr)\nonumber\\
   &+  \|u(t) \|_{s-4} + \| \pa u(t) \|_{s-4} 
   + \sup_{x \in \R^2} 
     \bigl\{ \langle t+|x| \rangle |u(t,x)|_{s-12}  \bigr\}\ 
 \Bigr]
\end{align*}
where $s\geq 29$ and $0<\delta <1$. Then we have the following. 

\begin{prp} \label{prp_apriori}
Assume that $f_j$, $g_j$ satisfy \eqref{DataSize}. Suppose that $E(T)\le 1$. 
There exists a positive constant $C_0$, which is independent of $\eps$ 
and $T$, such that 
\begin{align}
  E(T) \leq C_0(\eps + E(T)^2). 
  \label{mainest}
\end{align}
\end{prp}

\noindent{\bf Proof.}\ \ The following argument is almost the same as that of 
the previous works (\cite{katayama}, \cite{ott}, \cite{suna1}, etc.). 
First we note that Corollary \ref{cor_normalform} and the commutation relation 
(\ref{commute}) imply 
\begin{align}
(\square +m_j^2) Z^{\alpha}( u_j-\Lambda_j) = Z^{\alpha} (N_j + R_j) 
 \label{rewrite}
\end{align}
with
\begin{align}
 |Z^{\alpha}\Lambda_j(t,x)| 
 &\lesssim  
 |u|_{[|\alpha|/2]+3} (|u|_{|\alpha|+2}+|\pa u|_{|\alpha|+2}),
 \label{lambda_est}\\
 |Z^{\alpha} R_j(t,x)|
   &\lesssim  | u |_{[|\alpha|/2]+5}^2 
              \bigl( |u|_{|\alpha|+4} + |\pa u|_{|\alpha|+4} \bigr),
 \label{r_est}
\end{align}
and
\begin{align}
 |Z^{\alpha} N_j(t,x)|
  \lesssim 
     \frac{1}{\langle t+|x| \rangle}\bigl|u \bigr|_{[|\alpha|/2]+4} 
       \bigl( |u|_{|\alpha|+3} + |\pa u|_{|\alpha|+3} \bigr)
 \label{n_est}
\end{align}
by Lemma \ref{lem_null}. 
We use Lemma \ref{lem_decay} for (\ref{rewrite}) with $|\alpha| \leq s-12$ 
to obtain 
\begin{align*}
  \sum_{|\alpha|\le s-12}
  \langle t+|x|\rangle &|Z^{\alpha}(u_j -\Lambda_j)(t,x)|\\
  \lesssim  
   \eps &+
   \sum_{j=0}^{\infty} \sum_{|\beta|\leq s-8}
      \sup_{ \tau \in [0,t]} \varphi_j(\tau) 
      \bigl\| 
       \langle \tau +|\cdot| \rangle 
       Z^{\beta} \bigl( N_j+ R_j\bigr)(\tau,\cdot)
      \bigr\| _{L^2},
\end{align*}
where we have used
$$
\sum_{j=0}^\infty
\sum_{|\beta|\le s-7}
 \bigl\|
   \langle y\rangle \varphi_j(|y|)\bigl(Z^\beta(u_j-\Lambda_j)\bigr)(0,y)
 \bigr\|_{L^2(\R_y^2)}
 \lesssim 
 \sum_{j=0}^\infty 2^{-j}\bigl(\|f\|_{H^{s-4,s-2}}+\|g\|_{H^{s-5, s-2}}\bigr)
\lesssim \eps
$$
by \eqref{DataSize}. From (\ref{r_est}) and (\ref{n_est}) it follows that 
\begin{align*}
      \sum_{|\beta|\le s-8} &\bigl\| 
        \langle t + |\cdot|\rangle Z^{\beta} 
        (N_j+ R_j)(t,\cdot) \bigr\| _{L^2}\\
    &\lesssim  
       \frac{1}{\langle t \rangle} 
        \bigl( \|u(t)\|_{s-4} + \|\pa u(t)\|_{s-4} \bigr)
       \sup_{y \in \R^2} 
       \Bigl(
        \langle t +|y| \rangle \bigl| u(t,y) \bigr|_{[(s-8)/2]+5} 
       \Bigr) \\
 &\quad + 
       \frac{1}{\langle t \rangle} 
       \bigl( \|u(t)\|_{s-4} + \|\pa u(t)\|_{s-4} \bigr)
       \sup_{y \in \R^2} 
       \Bigl( 
        \langle t +|y| \rangle \bigl| u(t,y)\bigr|_{[(s-8)/2]+5} 
       \Bigr)^2 \\
 &\lesssim  
       \frac{E(T)^2}{\langle t\rangle}.
\end{align*}
Here we have used the relation $[(s-8)/2]+5 \leq s-12 $ for $s\geq 25$. 
So we have 
\begin{align*}
  \langle t+|x| \rangle |u_j(t,x) -\Lambda_j(t,x)|_{s-12} 
  &\lesssim
     \eps + \sum_{j=0}^{\infty} \sup_{\tau \in [0,t]} 
     \varphi_j(\tau) \frac{E(T)^2}{\langle \tau \rangle}\\
  &\lesssim
    \eps +E(T)^2 \sum_{j=0}^{\infty} 2^{-j}\\
  &\lesssim  
    \eps +E(T)^2.
\end{align*}
Also, (\ref{lambda_est}) and the Sobolev embedding theorem yield 
\begin{align*}
 \langle t+ |x| \rangle|\Lambda_j(t,x)|_{s-12} 
 &\lesssim  \langle t+ |x| \rangle|u(t,x)|_{[(s-12)/2]+3} |u(t,x)|_{s-9}\\
 &\lesssim \langle t+ |x| \rangle|u(t,x)|_{s-12} \|u(t)\|_{s-7}\\
 &\lesssim E(T)^2.
\end{align*}
Summing up, we obtain 
\begin{align}
 \langle t+|x| \rangle |u(t,x)|_{s-12} \lesssim \eps +E(T)^2 
  \label{pointwise}
\end{align}
for $(t,x) \in [0,T)\times \R^2$. 
Next we apply the standard energy inequality to (\ref{rewrite}) with 
$|\alpha|\le s-4$. Then we obtain 
\begin{align*}
 \| (u_j - \Lambda_j )(t) \|_{s-4} + 
 \|\pa (u_j - \Lambda_j )(t)\|_{s-4} 
 &\lesssim
  \eps + \int_{0}^{t} 
  \bigl\| (N_j + R_j)(\tau)\bigr\|_{s-4} d\tau .    
\end{align*}
By using  (\ref{r_est}) and (\ref{n_est}) again, we see that 
\begin{align*}
  \bigl\| (N_j + R_j)(t)\bigr\|_{s-4}
  &\lesssim  
    \frac{1}{ \langle t \rangle^2} 
    \bigl( \|u(t)\|_{s} + \|\pa u(t)\|_{s} \bigr)
       \sup_{y \in \R^2} \Bigl(
        \langle t +|y| \rangle  \bigl| u(t,y) \bigr|_{[\frac{s-4}{2}]+5} 
       \Bigr) \\
 &\quad + 
       \frac{1}{\langle t \rangle^{2}} 
       \bigl( \|u(t)\|_{s} + \|\pa u(t)\|_{s} \bigr)
       \sup_{y \in \times \R^2} 
       \Bigl( 
        \langle t +|y| \rangle \bigl| u(t,y)\bigr|_{[\frac{s-4}{2}]+5} 
       \Bigr)^2 \\
  &\lesssim  
       \frac{E(T)^2}{\langle t \rangle^{2-\delta}}. 
\end{align*}
Here we have used the relation $[(s-4)/2]+5 \le s-12$ for $s\geq 29$. 
So we have 
\begin{align*}
 \| (u_j - \Lambda_j )(t) \|_{s-4} + 
 \|\pa (u_j - \Lambda_j )(t)\|_{s-4} 
 \lesssim 
   \eps + \int_{0}^{t}\frac{E(T)^2}{\langle \tau \rangle^{2-\delta}}  d\tau
 \lesssim
   \eps +E(T)^2.    
\end{align*}
Also, (\ref{lambda_est}) leads to 
\begin{align*}
  \|\Lambda_j(t) \|_{s-3} + \| \pa \Lambda_j(t) \|_{s-4}
 &\lesssim 
 \frac{1}{\langle  t\rangle} 
 \sup_{y \in \R^2} \Bigl(\langle  t + |y| \rangle |u(t,y)|_{[(s-3)/2]+4}\Bigr)
 \bigl( \|u(t)\|_{s-1} + \|\pa u(t) \|_{s-1}\bigr)\nonumber \\
 &\lesssim 
  E(T)^2\langle t \rangle^{\delta-1}. 
\end{align*}
To sum up, we have
\begin{align}
 \|u(t)\|_{s-4} + \| \pa u(t) \|_{s-4}
 \lesssim \eps + E(T)^2 
 \label{low_energy}
\end{align}
for $t \in [0,T)$. 
Finally we apply $Z^{\alpha}$ to (\ref{eq}) with $|\alpha| \le s$ 
to obtain
\begin{equation}
(\Box+m_j^2)(Z^\alpha u_j)
-\sum_{k=1}^2\sum_{a,b=0}^3 \gamma_{ab}^{jk}(u,\pa u) \pa_a\pa_b (Z^\alpha u_k)
=F_j^{(\alpha)}, \quad j=1,2
\end{equation}
with 
$F_j^{(\alpha)}=Z^\alpha\bigl(F_j(u,\pa u, \pa^2 u)\bigr)
-\sum_{k=1}^2\sum_{a,b=0}^3 \gamma_{ab}^{jk}(u,\pa u) 
\pa_a\pa_b (Z^\alpha u_k)$, 
where $\gamma=(\gamma_{ab}^{jk})$ is from \eqref{QL}. 
Because of \eqref{Symmetry}, we can use the energy inequality for
hyperbolic systems with symmetric variable coefficients to estimate 
$\|Z^\alpha u(t)\|_{L^2}+\|\pa Z^\alpha u(t)\|_{L^2}$, and we see that 
\begin{align*}
 \|u(t) \|_{s}+ \|\pa u(t)\|_{s} 
  &\lesssim 
   \eps 
   + \int_{0}^{t} \bigl\|
       \pa\bigl(\gamma(u,\pa u)\bigr)(\tau)
     \bigr\|_{L^\infty}\|\pa u(\tau)\|_{s}
   +\|F^{(\alpha)}(\tau)\|_{s} d\tau.
\end{align*}
Since 
$$
 |F^{(\alpha)}| 
 \lesssim  |u|_{[s/2]+2} (|u|_{s}+|\pa u|_{s})
 \lesssim  |u|_{s-12} (|u|_{s}+|\pa u|_{s})
$$
and
$$
\bigl|\pa\bigl(\gamma_{ab}^{jk}(u,\pa u)\bigr)\bigr|\lesssim |u|_1+|\pa u|_1,
$$
we have 
\begin{align}
 \|u(t) \|_{s}+ \|\pa u(t)\|_{s} 
  &\lesssim 
   \eps + \int_{0}^{t} E(T)^2 \langle \tau \rangle^{\delta -1}d\tau 
   \nonumber \\
 &\lesssim  
   \eps + E(T)^2 \langle t \rangle^{\delta} 
 \label{high_energy}
\end{align}
for $t \in [0,T)$. 
By (\ref{pointwise}), (\ref{low_energy}) and (\ref{high_energy}), 
we arrive at the desired estimate (\ref{mainest}).
\qed\\

Now we finish the proof of Theorem \ref{thm_main}. 
The inequality (\ref{mainest}) implies 
that there exists a constant $M>0$, which does not depend on $T$, such that 
$$E(T) \le M$$ 
if we choose $\eps$ sufficiently small. 
The unique global existence is an immediate consequence of 
this {\it a priori} bound and the classical local existence theorem 
(see \cite{hor} etc.). 
To prove the existence of a free profile, we remember that 
$$
 (\Box+m_j^2)(u_j-\Lambda_j) = N_j +R_j
$$
with 
\begin{align*}
 &\|(N_j+R_j)(t,\cdot)\|_{H^{s-4}(\R^2)} 
 \lesssim \langle t \rangle^{-2+\delta},
 \\
 &\|\Lambda_j(t,\cdot)\|_{H^{s-3}(\R^2)} 
 + \|\pa_t \Lambda_j(t,\cdot) \|_{H^{s-4}(\R^2)}
 \lesssim \langle t \rangle^{-1+\delta}.
\end{align*}
Now we set
\begin{align*} 
 &f_j^+
 =
 f_j-\Lambda_j\bigr|_{t=0} 
 +
 \int_{0}^{\infty} 
 \frac{\sin \left(-\tau \OOmega_j \right)}{\OOmega_j }
 (N_j+R_j)(\tau,\cdot) d\tau,\\
 &g_j^+
 =
 g_j-\pa_t \Lambda_j \bigr|_{t=0} 
 +
 \int_{0}^{\infty} 
 \left(\cos(-\tau \OOmega_j )\right) (N_j+R_j)(\tau,\cdot) d\tau
\end{align*}
and
$$
 u_j^+(t,\cdot)= 
 \left(\cos(t\OOmega_j )\right) f_j^+
 +
 \frac{\sin \left(t \OOmega_j \right)}{\OOmega_j }g^+_j
$$
with $\OOmega_j=(m_j^2-\Delta)^{1/2}$. 
Since the Duhamel formula yields
\begin{align*}
 u_j(t,\cdot)-\Lambda_j(t,\cdot) 
 =&
 \left(\cos(t\OOmega_j )\right) (f_j - \Lambda_j|_{t=0})
 +
 \frac{\sin \left(t \OOmega_j \right)}{\OOmega_j }
 (g_j -\pa_t \Lambda_j|_{t=0})\\
 &+
  \int_{0}^{t} 
 \frac{\sin \left((t-\tau) \OOmega_j  \right)}{\OOmega_j }
 (N_j+R_j)(\tau,\cdot) d\tau
 \\
 =& 
 u_j^{+}(t,\cdot)-
  \int_{t}^{\infty} 
 \frac{\sin \left((t-\tau) \OOmega_j  \right)}{\OOmega_j }
 (N_j+R_j)(\tau,\cdot) d\tau,
\end{align*}
we have
\begin{align*}
 &\|(u_j-u_j^+)(t,\cdot)\|_{H^{s-3}(\R^2)}
 +\|\pa_t (u_j-u_j^+)(t,\cdot)\|_{H^{s-4}(\R^2)}\\
 &\lesssim
 \|\Lambda_j(t,\cdot)\|_{H^{s-3}(\R^2)}
 +
 \|\pa_t \Lambda_j(t,\cdot)\|_{H^{s-4}(\R^2)}
 +
 \int_{t}^{\infty} \|(N_j+R_j)(\tau,\cdot) \|_{H^{s-4}(\R^2)}d\tau\\
 &\lesssim
 \langle t \rangle^{-1+\delta} 
 + 
 \int_{t}^{\infty} \langle \tau \rangle^{-2+\delta}d\tau\\
 &\lesssim
  \langle t \rangle^{-1+\delta}.
\end{align*}
This completes the proof of Theorem \ref{thm_main}.\qed\\

\medskip
\subsection*{ Acknowledgments}

The authors would like to express their sincere gratitude to Professor 
Jalal Shatah for the fruitful discussion that motivates the present work,
and also for his warm hospitality during their visit to the Courant Institute 
of Mathematical Sciences, New York University, 
where a part of this work was done.

The first author (S.K.) is partially supported by Grant-in-Aid for Scientific 
Research (C) (No.20540211), JSPS. 
The second author (T.O.) is partially supported by Grant-in-Aid for Scientific 
Research (A) (No.21244010), JSPS. 
The third author (H.S.) is partially supported by Grant-in-Aid for Young 
Scientists (B) (No.22740089), MEXT.

\end{document}